\documentclass[11pt]{article}

\usepackage{graphicx}
\usepackage{amsmath}
\usepackage{amsfonts}
\usepackage{fancyhdr} 

\usepackage{wrapfig}

\usepackage{bbold}

\usepackage{float}

\newtheorem{lem}{Lemma}
\newtheorem{thm}{Theorem}

\newtheorem{prop}[lem]{Proposition}

\numberwithin{equation}{section}
\numberwithin{thm}{section}
\numberwithin{lem}{section}
\numberwithin{defn}{section}
\numberwithin{alg}{section}

\newcommand{\sss}{\mathsf{S}}

\newcommand{\real}{\mathbb{R}}
\newcommand{\integer}{\mathbb{Z}}

\newenvironment{pf}{\noindent {\em Proof}.\ \ }{\hspace*{\fill}\rule{.5ex}{1.4ex}\,}

\newcommand{\nat}{\mathbb{N}}

\newcommand{\gtr}{G^{(\tau,R)}}

\newcommand{\mpp}{\mathsf{p}}

\title{The combinatorics of scattering in layered media}
\author{Peter C.~Gibson \footnote{Dept.~of Mathematics \& Statistics, York University, 4700 Keele St., Toronto, Ontario, Canada, M3J~1P3, $\mathtt{pcgibson@yorku.ca}$} 
}

\begin{document}

\date{May 17, 2013}
\maketitle

\begin{abstract} 
Reflection and transmission of waves in piecewise constant layered media are important in various imaging modalities and have been studied extensively.  Despite this, no exact time domain formulas for the Green's functions have been established.  Indeed, there is an underlying combinatorial obstacle: the analysis of scattering sequences.   In the present paper we exploit a representation of scattering sequences in terms of trees to solve completely the inherent combinatorial problem, and thereby derive new, explicit formulas for the reflection and transmission Green's functions.  
\end{abstract}

%
%

\section{Introduction}
The analysis of wave propagation in layered media is a well-developed subject, with applications to acoustic, electromagnetic and seismic imaging; see \cite{BrGo:1990}, \cite{FoGaPaSo:2007}, \cite{BlCoSt:2001} and the numerous references therein.  In the present introduction we describe briefly the physical framework and cite some needed facts. In the next section we explain a fundamental combinatorial problem arising from the given framework---the enumeration of scattering sequences.  The main contribution of the present paper is to solve the combinatorial problem completely, thereby producing new, exact formulas for the reflection and transmission Green's functions. Our main results are Theorem~\ref{thm-reflection} and Theorem~\ref{thm-transmission}. 

\begin{figure}[h]
\parbox{3.3in}{
\fbox{
\includegraphics[clip,trim=.75in 0in 0in 0in, width=3in]{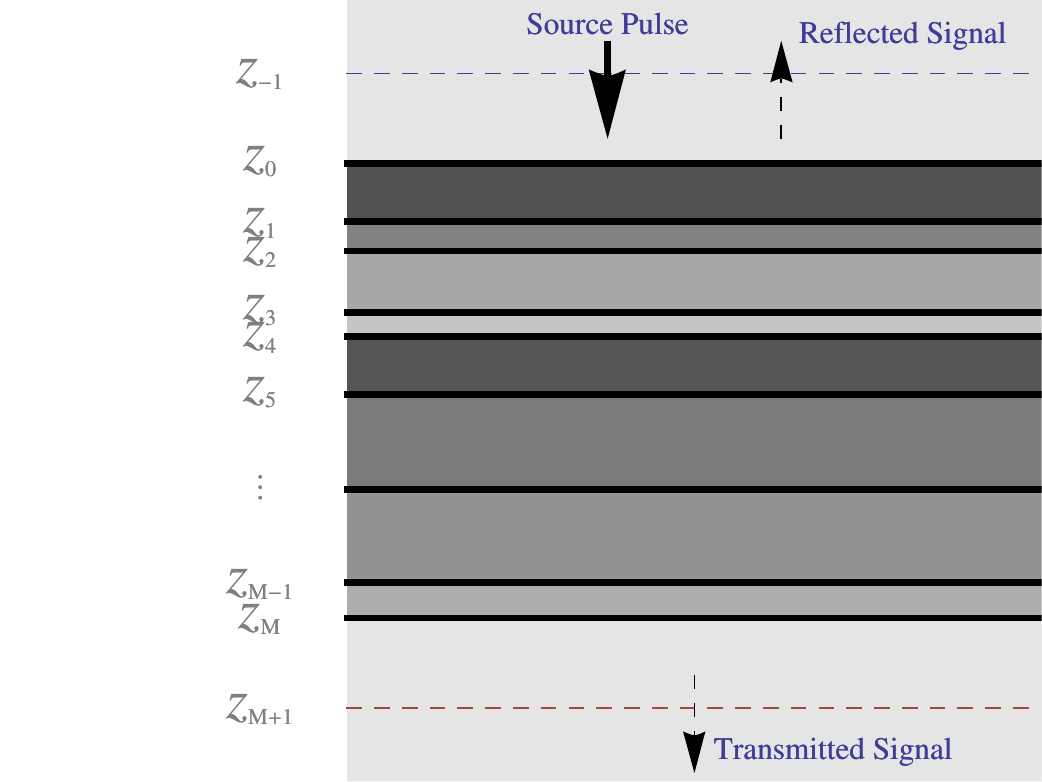}
}
\caption{ A layered medium with $M$ layers, sitting between two half spaces. The depth coordinate $z$ increases downward.}\label{fig-Layered}
}
\end{figure}

The basic framework is as follows.  Let $(x,y,z)$ denote euclidean coordinates for a solid three-dimensional acoustic medium whose physical parameters (density and bulk modulus) vary only in the $z$-direction.   Suppose furthermore that these physical parameters are piecewise constant in $z$, having jumps at the $M+1$ locations
\[
z_0<z_1<\cdots<z_M.
\]
Thus the solid contains $M$ homogenous layers, the $n$th layer corresponding to the $z$-interval
\[
z_{n-1}<z<z_n;
\]
the layers are sandwiched between the two half spaces $z<z_0$ and $z>z_M$.   We shall refer to the $z$-direction as depth, and depict it as increasing downward, as in Figure~\ref{fig-Layered}.  

\pagestyle{fancyplain}

From the perspective of acoustic imaging the problem is to infer the physical parameters for the $M$ layers by probing them with acoustic waves in one of two ways.   The first way is to send an acoustic pulse from a fixed depth $z_{-1}<z_0$ toward the interface at $z_0$ and to record the pulse train reflected back from the $M$ layers as it crosses the original depth $z_{-1}$.   The \emph{reflection problem} is to infer the structure of the $M$ layers from this recorded reflection data.  A second experiment is to transmit a pulse from $z_{-1}$ as before, but to record the resulting pulse train that is transmitted across the $M$ layers to some fixed depth $z_{M+1}>z_M$.   The \emph{transmission problem} is to infer the structure of the $M$ layers from this transmitted data, recorded at $z_{M+1}$.   It is assumed for both the reflection and transmission experiments that the initial pulse is a plane wave of the form $f(z-ct)$, where $t$ denotes time and $c(z)$ denotes the speed of sound at depth $z$.  Thus although the physical setting is three-dimensional, by symmetry the analysis reduces to one spatial dimension, the $z$-direction.  

If the initial pulse is idealized to a Dirac delta function initially centred at $z_{-1}$, so of the form
\[
\delta(z-z_{-1}-ct),
\]
then the reflection data $G(t)$ and the transmission data $H(t)$ are Green's functions for the governing partial differential equation (see Appendix~\ref{sec-appendix}), and they have the form
\begin{align}
G(t)=& \sum_{n=1}^\infty a_n\,\delta(t-\sigma_n)\label{G},\\
H(t)\rule{0pt}{24pt}=& \sum_{n=1}^\infty b_n\,\delta(t-\sigma_n)\label{H}.
\end{align}
The numbers $a_n$ and $b_n$ will be referred to as reflection amplitudes and transmission amplitudes, respectively; the numbers $\sigma_n$ and $\sigma_n^\prime$ will be referred to as arrival times.   

As a consequence of the standard theory, both $G$ and $H$ are completely determined by two vectors:
\begin{enumerate}
\item the sequence $R=(R_0,\ldots,R_M)$ of reflection coefficients at the layer boundaries; and
\item the sequence $\tau^\prime=(\tau_0,\ldots,\tau_{M+1})$ of two-way travel times, where $\tau_n$ denotes twice the time it takes a downward-traveling acoustic wave to go from $z_{n-1}$ to $z_n$.   
\end{enumerate}
The travel time $\tau_{M+1}$ is obviously irrelevant for the reflection problem, so $G$ is in fact determined by $R$ and $\tau=(\tau_0,\ldots,\tau_M)$.   We incorporate this basic fact into our notation, writing
\[
G^{(\tau,R)}\quad\mbox{and}\quad H^{(\tau^\prime,R)}
\]
to denote the respective reflection and transmission data determined by a given pair $(\tau^\prime,R)$.    

Different layered media that correspond to the same parameters $(\tau^\prime,R)$ are for present purposes indistinguishable, and we shall simply refer to the pair $(\tau^\prime,R)$ as a medium, letting it be understood that a class of physical systems are thereby represented.  (See Appendix~\ref{sec-appendix}.)

The above facts lead naturally to a basic question: What is the formula for the amplitude coefficients $a_n$ and $b_n$ in terms of $(\tau^\prime,R)$?   It turns out that no finite, closed form formula has ever been established---only series expansions that must be estimated when it comes to practical computation.   (See, for example, \cite[Chapter~3]{FoGaPaSo:2007}, \cite[Chapter~2]{BlCoSt:2001}, \cite[Section~2.5]{BrGo:1990}; concerning theoretical developments applied specifically to seismic, see \cite{Ne:1980}, \cite{SaSy:1988}, \cite{WeArFe:2003},\cite{In:2009})    This is because there is a substantial combinatorial problem blocking the way to an exact formula.   The purpose of the present paper is to solve the combinatorial problem directly, and to present the resulting explicit formulas for the reflection and transmission Green's functions $G$ and $H$ in terms of $(\tau^\prime,R)$.

\section{The combinatorics of reflection\label{sec-reflection}}

The one-dimensional wave equation for a homogenous medium generates traveling wave solutions; in a layered medium one has to take into account the behaviour of these traveling waves at interfaces between homogenous layers.  This behaviour depends on the reflection coefficient $R_n$ associated with the interface at $z_n$ in the following way.  When a wave $f(z-z_n-c_nt)$ traveling downward from $z_{n-1}$ toward $z_n$ for $t<0$ hits the interface at $t=0$, it splits into a reflected wave,
\[
R_nf(z-z_n+c_nt)
\]
that travels back up toward $z_{n-1}$, and a transmitted wave
\[
\sqrt{1-R_n^2}\,f(z-z_n-c_{n+1}t)
\] 
that continues down toward $z_{n+1}$ (at modified speed $c_{n+1}$).   The reflected and transmitted waves then hit the respective interfaces at $z_{n-1}$ and $z_{n+1}$, generating two new sets of reflected and transmitted waves, and a cascade of successive reflections and transmissions continues indefinitely.  The case of an initial wave traveling upward from $z_{n+1}$ toward $z_n$ instead of downward from $z_{n-1}$ is similar, except that the reflection coefficient $-R_n$ applies instead of $R_n$; the transmission coefficient $\sqrt{1-R_n^2}$ remains unchanged.   We use the notation $T=(T_0,\ldots,T_M)$ for the transmission coefficients that correspond to reflection coefficients $R=(R_0,\ldots,R_M)$ by way of the formula
\[
T_n=\sqrt{1-R_n^2}\quad\quad(0\leq n\leq M).
\]

A standard idea is to view $\gtr(t)$ as the sum total of all possible sequences of successive reflections and transmissions of an initial pulse at $z_{-1}$ that eventually return to $z_{-1}$.  (See \cite[Chapter~3]{FoGaPaSo:2007} for details.) 
The idea is worth illustrating further, since it underpins the arguments below. For example, consider an initial downward traveling unit pulse of the form $\delta(z-z_{-1}-c_0t)$.  After $\tau_0/2$ seconds, the initial pulse  
reaches the interface at $z_0$.   Part of it is transmitted into the first layer as $T_0\delta\bigl(z-z_0-c_1(t-\frac{\tau_0}{2})\bigr)$, which, after another $\tau_1/2$ seconds, reaches the interface $z_1$.   Part of this pulse is reflected back into first layer as $R_1T_0\delta\bigl(z-z_1+c_1(t-\frac{\tau_0+\tau_1}{2})\bigr)$.   Traveling back up to $z_0$, the latter reaches $z_0$ after a further $\tau_1/2$ seconds, and is partly transmitted back to $z_{-1}$ as 
\begin{equation}\label{samplewave}
R_1T_0^2\delta\bigl(z-z_0+c_0(t-\textstyle\frac{\tau_0+2\tau_1}{2})\bigr),
\end{equation}
arriving at $z_{-1}$ at time $\sigma=\tau_0+\tau_1$.   Thus the part of the initial pulse that traverses the sequence of depths $\mpp=(z_{-1},z_0,z_1,z_0,z_{-1})$ returns to $z_{-1}$ with an amplitude $\alpha=R_1T_0^2$ at arrival time $\sigma=\tau_0+\tau_1$, thereby contributing a term of the form $\alpha\delta(t-\sigma)$ to $\gtr(t)$.   The sequence $\mpp$ is called a scattering sequence, and the amplitude $\alpha$ is called the weight of $\mpp$.   As mentioned earlier, the impulse response itself is a delta train of the form
\begin{equation}\label{normalform}
G^{(\tau,R)}(t)=\sum_{n=1}^\infty\alpha_n\delta(t-\sigma_n),
\end{equation}
composed of the cumulative contributions of all possible scattering sequences returning to $z_{-1}$, with their associated weights and arrival times.  In the above example, $\mpp$ is the only scattering sequence having arrival time $\sigma$.  But in general different scattering sequences may arrive simultaneously, so each amplitude $\alpha_n$ occurring in (\ref{normalform}) is the sum of the weights of all scattering sequences arriving at time $\sigma_n$.   Given an arrival time $\sigma_n$, the essential combinatorial problem, then, is to enumerate all the associated scattering sequences, along with their weights, to derive a formula for $\alpha_n$.

\subsection{Scattering sequences and the Green's function}

The first step is to establish some terminology and notation, as follows.  A scattering sequence that starts and ends at $z_{-1}$ is represented by a path in the graph
\begin{equation}\label{simplegraph}
\mbox{
\includegraphics[width=4in]{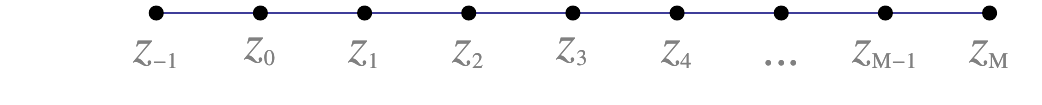}
}
\end{equation}
Formally, given an integer $M\geq 1$, let $\sss_M$ denote the set of all sequences of the form
\[
\mpp=(\mpp_0,\mpp_1,\ldots,\mpp_L)
\]
such that $L\geq 2$ and: 
\begin{subequations}
\begin{align}
\mpp_0=\mpp_L=z_{-1}\mbox{ and }\mpp_n\in&\left\{z_1,\ldots,z_M\right\}\mbox{ if }1\leq n\leq L-1;\label{pathcondition1}\\
\begin{split}
\forall n \mbox{ with } 0\leq n\leq L-1,&\quad\exists j \mbox{ with }-1\leq j\leq M-1\mbox{ such that }\\
&\{\mpp_n,\mpp_{n+1}\}=\{z_j,z_{j+1}\}.\label{pathcondition2}
\end{split}
\end{align}
\end{subequations}
The elements of $\sss_M$ will be referred to as scattering sequences.  Condition (\ref{pathcondition1}) says a scattering sequence starts and ends at $z_{-1}$, and condition (\ref{pathcondition2}) (which refers to \emph{un}ordered pairs) says that adjacent terms in a scattering sequence are adjacent vertices in the graph (\ref{simplegraph}). 

For example, for $0\leq n\leq M$, the shortest scattering sequence that reaches $z_n$ is  
\begin{equation}\label{primary}
(z_{-1},z_0,z_1,\ldots,z_{n-1},z_n,z_{n-1},\ldots,z_1,z_0,z_{-1});
\end{equation} 
this is called a primary scattering sequence.  A scattering sequence reaching maximum depth $z_n$ that is not shortest possible is called a multiple scattering sequence.

\subsubsection{The weight of a scattering sequence\label{sec-weight}}

The weight $w(\mpp)$ corresponding to a scattering sequence 
\[
\mpp=(\mpp_0,\mpp_1,\ldots,\mpp_L)
\]
in an $M$-layer medium $(\tau,R)$ is defined as follows.   For each $n$ in the range $1\leq n\leq L-1$, and given that $\mpp_n=z_j$, define
 \begin{equation}\label{wn}
 w_n=\left\{\begin{array}{cc}
 R_j&\mbox{ if }\mpp_{n-1}=\mpp_{n+1}=z_{j-1}\\
 -R_j&\mbox{ if }\mpp_{n-1}=\mpp_{n+1}=z_{j+1}\\
 \sqrt{1-R_j^2}&\mbox{ otherwise }
 \end{array}\right..
 \end{equation}
 The three possibilities correspond respectively to: reflection inside the $j$th layer at $z_j$; reflection inside the $(j+1)$st layer at $z_j$; and transmission between the $j$th and $(j+1)$st layers.  Finally, set 
 \begin{equation}\nonumber
 w(\mpp)=\prod_{n=1}^{L-1}w_n.
 \end{equation}
Thus the part of an initial unit impulse that traverses $\mpp$ returns to $z_{-1}$ with amplitude $w(\mpp)$.

\subsection{Transit count and branch count vectors}

We define two maps, 
\[
\kappa,\beta:\sss_M\rightarrow\integer^{M+1}
\]
that associate integer vectors to a given scattering sequence.    

\begin{figure}[h]
\fbox{
\includegraphics[clip,trim=0in 1in 0in 0in, width=4.8in]{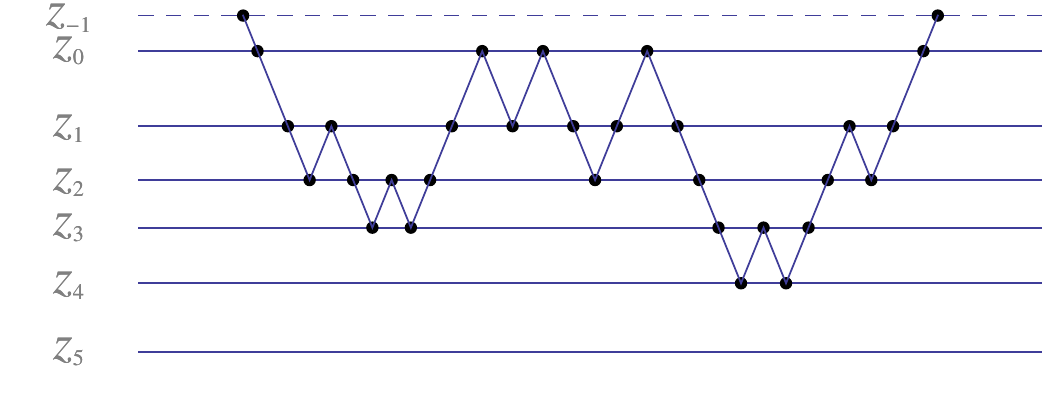}
}
\caption{ The Dyck path for a scattering sequence 
$
\mpp,
$
with the horizontal coordinate now representing time, which increases to the right.}\label{fig-Dyck}
\end{figure}

A scattering sequence $\mpp\in\sss_M$  may be represented graphically as in Figure~\ref{fig-Dyck};  Stanley \cite{St:1999}  calls  such a representation a Dyck path.   
\begin{figure}[h]
\fbox{
\includegraphics[clip,trim=0in .9in 0in 0in, width=4.8in]{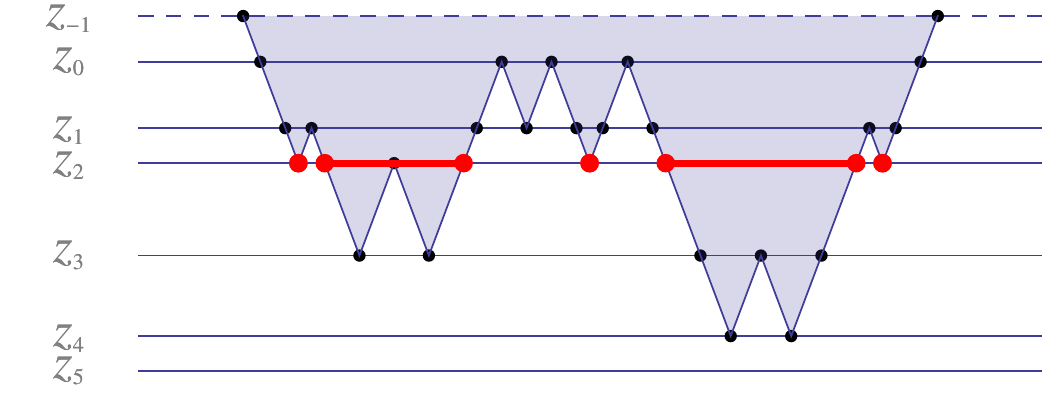}
}
\caption{ The intervals $I^2_j$, in red.  There are five intervals, two of which have positive length, so $k_2=5$ and $b_2=2$.}\label{fig-Intervals}
\end{figure}

Given the Dyck path for a scattering sequence $\mpp\in\sss_M$, let $t$ denote the horizontal coordinate and (as usual) let $z$ denote the vertical coordinate.   Let $U$ denote the portion of the $t,z$-plane on or above the Dyck path and at or below $z=z_{-1}$---the shaded region in Figure~\ref{fig-Intervals}.  For each $n$ in the range $0\leq n\leq M$, consider the horizontal line $L_n$ at depth $z_n$.  The intersection $L_n\cap U$ consists of a disjoint union of closed intervals $I^n_j$, where $1\leq j\leq k_n$; see Figure~\ref{fig-Intervals}.   The intervals of $I^n_j$ are of two types: degenerate intervals consisting of a single points; and non-degenerate intervals having positive length.    Letting $k_n$ denote the total number of intervals $I^n_j$, and letting $b_n$ denote the number of non-degenerate intervals, set 
\[
\kappa(\mpp)=k=(k_0,\ldots,k_M)\quad\mbox{ and }\quad\beta(\mpp)=b=(b_0,\ldots,b_M).
\]

Observe that the entry $k_n$ of the vector $k=\kappa(\mpp)$ counts the number of times the Dyck path crosses back and forth across the $n$th layer $z_{n-1}<z<z_n$; the vector $k$ is therefore called the transit count vector for $\mpp$.   The vector $b=\beta(\mpp)$ is called the branch count vector, for reasons that will be apparent in Section~\ref{sec-tree}.

\subsubsection{The set of transit count vectors}

The range of $\kappa$ is contained in the set
\begin{equation}\label{LM}
\begin{split}
&\mathfrak{L}_M=\\
&\left\{(k_0,k_1,\ldots,k_M)\in\nat^{M+1}\,\left|\, k_0=1\,\&\,\forall n\leq M-1,\,k_n=0\Rightarrow k_{n+1}=0\right.\right\}.
\end{split}
\end{equation} 
Conversely, for any given $k\in\mathfrak{L}_M$, it is straightforward to construct a realizing scattering sequence.  Thus $\mathfrak{L}_M$ is precisely the range $\kappa(\sss_M)$ of the mapping $\kappa:\sss_M\rightarrow\integer^{M+1}$.

\subsubsection{Arrival times}

The arrival time of $\mpp$, expressed in terms of the transit count vector $k=\kappa(\mpp)$, is simply 
\[
\langle k,\tau\rangle=k_0\tau_0+k_1\tau_1+\cdots+k_M\tau_M.
\]
Therefore $\gtr(t)$ may be written as 
\begin{equation}\label{arrival}
\gtr(t)=\sum_{k\in\mathfrak{L}_M}a(R,k)\delta(t-\langle k,\tau\rangle),
\end{equation}
where the amplitudes $a(R,k)$ are given by the formula
\begin{equation}\label{aRkdefintion}
a(R,k)=\sum_{\substack{\mpp\in\sss_M\; \scriptstyle{\rm such\; that}\\ \kappa(\mpp)=k}}w(\mpp).
\end{equation}
Note that the weight $w(\mpp)$ of a scattering sequence depends only on $R$, and not on $\tau$; the next step is to find an explicit formula for $a(R,k)$.   This is greatly facilitated by introducing another representation for scattering sequences, in terms of trees.  

\subsection{The tree representation of a scattering sequence\label{sec-tree}}

\begin{figure}[h]
\parbox{4.1in}{
\fbox{
\includegraphics[clip,trim=.9in .9in 0in 0in, width=3.8in, angle=180]{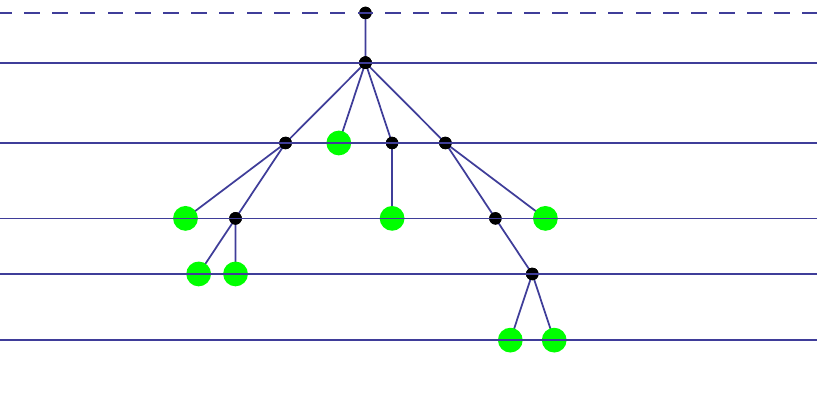}
}
\caption{A tree, with the leaves coloured green. The root is at the bottom, and horizontal lines indicate the various heights of vertices.}\label{fig-aTree}
}
\end{figure}

A tree is a connected cycle-free graph.  The vertices of a tree are divided into three anatomical types, as follows: the root is a single, specially designated vertex; a non-root that belongs to just one edge is called a leaf; all other non-root vertices are called branch points.     Vertices in a tree have a height, determined by their distance (in the sense of shortest path) to the root. See Figure~\ref{fig-aTree}. 

\begin{figure}[h]
\fbox{\parbox{\textwidth}{
(a)\includegraphics[clip,trim=0in .7in 0in 0in, width=2.25in]{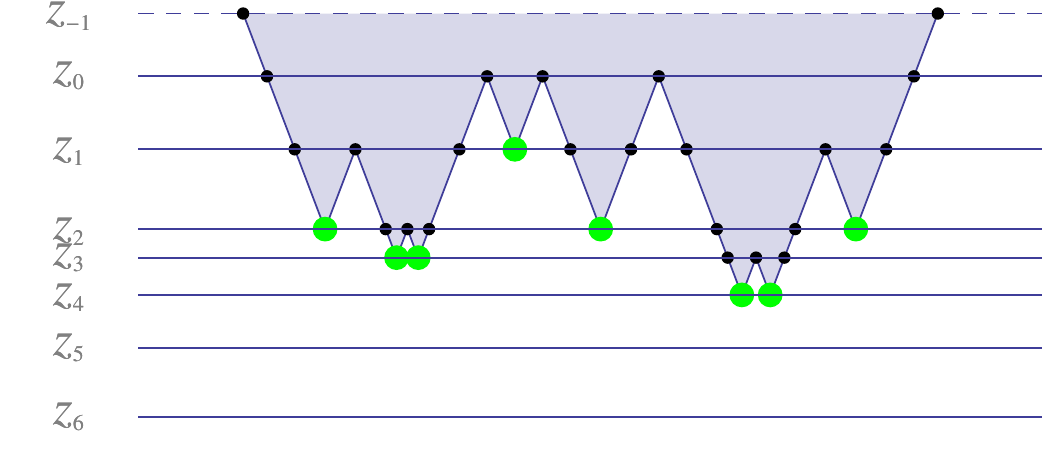}
(b)\includegraphics[clip,trim=0in .7in 0in 0in, width=2.25in]{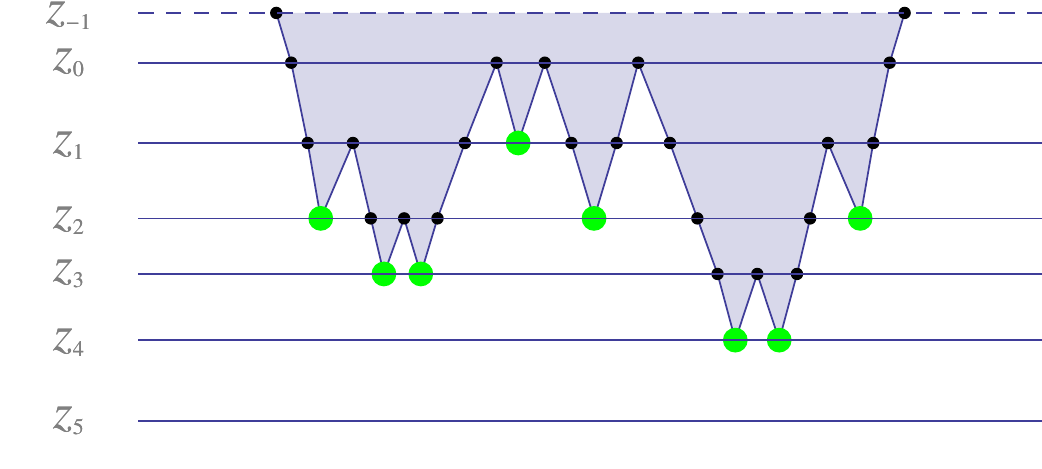}
(c)\includegraphics[clip,trim=0in .7in 0in 0in, width=2.25in]{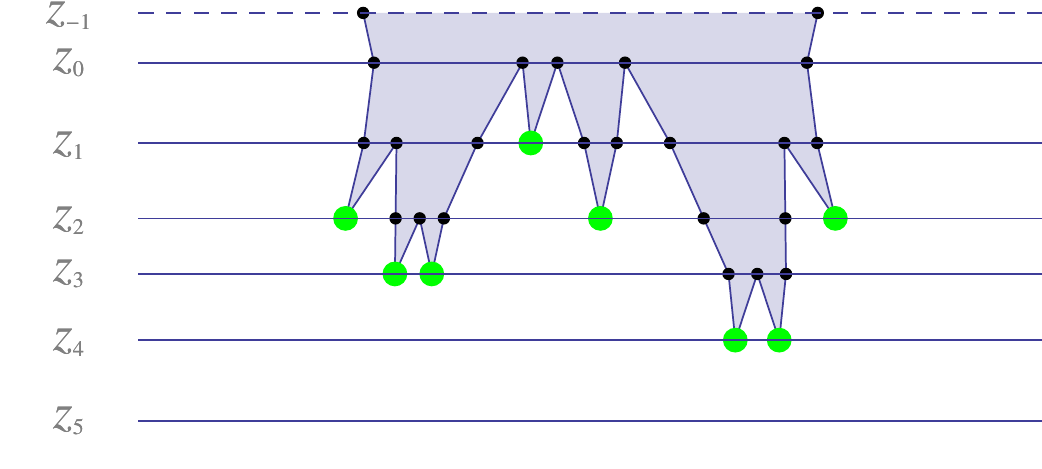}
(d)\includegraphics[clip,trim=0in 0.7in 0in 0in, width=2.25in]{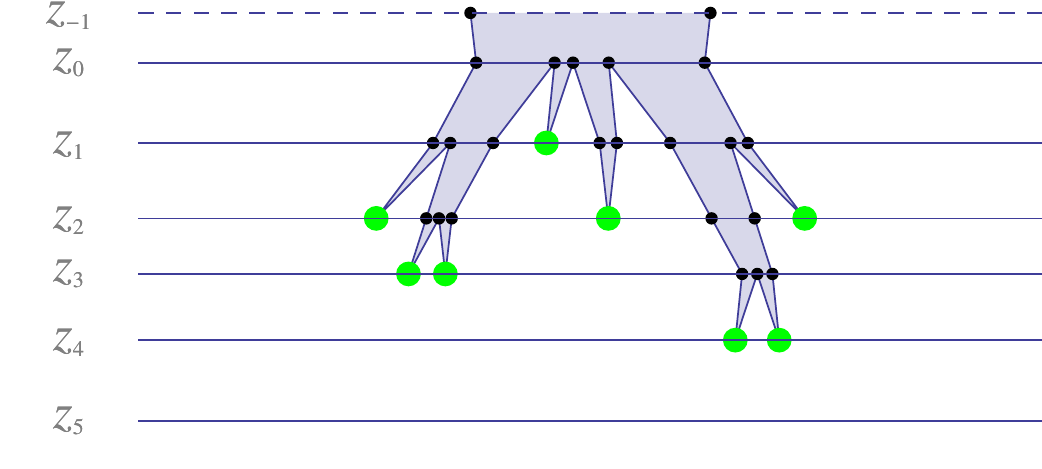}
(e)\includegraphics[clip,trim=0in 0.7in 0in 0in, width=2.25in]{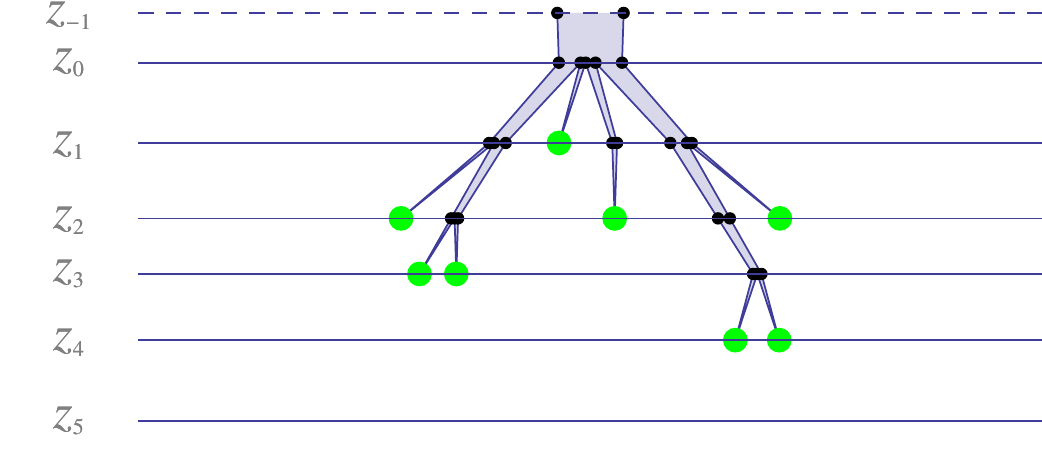}
(f)\includegraphics[clip,trim=0in 0.7in 0in 0in, width=2.25in]{Tree1.pdf}
}
}
\caption{ The collapsing of a Dyck path (a) to a tree (f) (depicted upside down), an operation which is reversible.}\label{fig-tree}
\end{figure}

The association between a scattering sequence and a tree is well-known (see \cite[Exercise 6.19]{St:1999}) and arises, for instance, in the analysis of  Brownian excursions and superprocesses (see \cite[Section~1.1]{Le:2005}).    The tree representing a scattering sequence $\mpp\in\sss_M$ may be obtained simply by collapsing its Dyck path, as follows.   Recall the intervals $I^n_j$ used to define $\kappa(\mpp)$ and $\beta(\mpp)$ above; for present purposes let $I^{-1}_1$ denote the intersection of $z=z_{-1}$ with $U$.  To collapse the Dyck path, contract each of the intervals $I^n_j$ to a point, keeping the distances between intervals unchanged, and interpolate this horizontal contraction linearly on each depth interval $z_{n-1}<z<z_{n}$.  This operation transforms the original Dyck path into a tree  (in fact it is an isotopy between the region $U$ and the resulting tree).  See Figure~\ref{fig-tree}.   Note that the degenerate intervals $I^n_j$ (coloured green for emphasis) end up as leaves, while non-degenerate intervals are contracted to branch points of the tree---except for $I^{-1}_1$, which is contracted to the root.    

Conversely, given a tree, one may recover the original scattering sequence by tracing the outline of the tree, from the root (keeping the tree on the left), and recording the depths of the vertices in the order that they are passed.  

Observe that the vectors $(k,b)=(\kappa(\mpp),\beta(\mpp))$ have a simple interpretation in terms of the tree representing $\mpp\in\sss_M$.   For each $0\leq n\leq M$, $k_n$ is the number of vertices at depth $z_n$, and $b_n$ is the number of branch points at $z_n$.  (This is the reason for calling $b$ the branch count vector for $\mpp$.)   There are some evident constraints on these quantities.  Note first that $b_M=0$.   Furthermore, for $0\leq n\leq M-1$, 
\begin{equation}\label{constraints}
\min\{1,k_{n+1}\}\leq b_n\leq k_{n+1}\quad(0\leq n\leq M-1),
\end{equation}
since each vertex at $z_{n+1}$ is connected by an edge to a unique branch point at $z_n$.    It is convenient to refer to the left shift $\tilde{k}$ of $k=\kappa(\mpp)$; that is,
\begin{equation}\label{ktilde}
\tilde{k}=(k_1,k_2,\ldots,k_M,0)\in\integer^{M+1}.
\end{equation}
In terms of this notation, the constraints
 (\ref{constraints}) become
\begin{equation}\label{constraints2}
\min\{1,\tilde{k}_n\}\leq b_n\leq\min\bigl\{k_n,\tilde{k}_n\bigr\}\quad(0\leq n\leq M).
\end{equation}

An easy application of the tree representation is to determine the possible values of $b=\beta(\mpp)$ for each given $k=\kappa(\mpp)$, or in other words, to determine the set 
\begin{equation}\label{V}
V(k)=\beta(\kappa^{-1}(k)).
\end{equation}
In fact $V(k)$ is determined precisely by (\ref{constraints2}). 
\begin{prop}\label{prop-V}
Given $k\in\mathfrak{L}_M$, 
\[
V(k)=\bigl\{b\,\bigl|\,\min\{\mathbb{1},\tilde{k}\}\leq b\leq\min\{k,\tilde{k}\}\bigr\};
\]
equivalently, $V(k)$ may be expressed as a Cartesian product of sets,
\[
V(k)=V_0\times V_1\times\cdots\times V_M,\]
where 
$
V_n=\bigl\{\min\{1,\tilde{k}_n\},\min\{1,\tilde{k}_n\}+1,\ldots,\min\{k_n,\tilde{k}_n\}\bigr\}\quad(0\leq n\leq M).
$
\end{prop}
Here $\mathbb{1}\in\integer^{M+1}$ is the vector whose entries are all 1.  The minimum is to be interpreted entrywise, meaning that for $x,y\in\real^{M+1}$, 
\[
\min\{x,y\}=\bigl(\min\{x_0,y_0\},\min\{x_1,y_1\},\ldots,\min\{x_M,y_M\}\bigr)\in\real^{M+1}.
\]
\begin{pf}
The constraints (\ref{constraints2}) imply that $V(k)\subset V_0\times\cdots\times V_M$.   It remains to show that any $b\in V_0\times\cdots\times V_M$ belongs to $V(k)$, which entails showing that there exists a scattering sequence $\mpp$ such that $(\kappa(\mpp),\beta(\mpp))=(k,b)$.   It suffices to construct a tree representing such a $\mpp$, as follows.  Given $b\in V_0\times\cdots\times V_M$, place $k_n$ vertices, consisting of $b_n$ branch points and $k_n-b_n$ leaves (in any order), at depth $z_n$, for $0\leq n\leq M$, and place a root at $z_{-1}$ (considered as a branch point).  Let $N$ denote the largest index such that $k_N\neq 0$.  For $0\leq n\leq N$,
\[
1\leq b_{n-1}\leq\min\{k_{n-1},k_n\},
\]
so there exists a surjection $f_n$ from the $k_n$ vertices at $z_n$ onto the $b_{n-1}$ branch points at $z_{n-1}$.   The tree representing $\mpp$ is completed by drawing an edge from each vertex $v$ at $z_n$ to $f_n(v)$.   
\end{pf}

\subsubsection{A formula for the weight}

The tree representation facilitates deriving a simple formula for the weights.  The following lemma uses multi-index notation, whereby given a vector $s=(s_0,s_1,\ldots,s_M)$ and an integer vector $d=(d_0,\ldots,d_M)$,  
\[s^d=\prod_{n=0}^Ms_n^{d_n}.\]  

\begin{lem}\label{lem-weight}
Let $\mpp\in\sss_M$ be a scattering sequence in an $M$-layer medium $(\tau,R)$, and set $(k,b)=(\kappa(\mpp),\beta(\mpp))$.  Then
\[
w(\mpp)=(-R)^{\tilde{k}-b}R^{k-b}T^{2b}.
\]   
\end{lem}
\begin{pf}
Consider the tree representing $\mpp$.  Observe that each instance of $w_n=R_j$ in (\ref{wn}) corresponds to a unique leaf at $z_j$ (see Figure~\ref{fig-tree}).  Since there are $k_n-b_n$ leaves at $z_n$ this results in a total contribution of $R^{k-b}$.   

Let $v$ be a branch point at $z_j$ having $d_v$ edges to vertices at $z_{j+1}$.   Observe that precisely $d_v-1$ of these edges correspond to an instance of $w_n=-R_j$ in (\ref{wn}), and every occurrence of $w_n=-R_j$ arises this way.  The sum total of numbers $d_v-1$ over branch points $v$ at $z_j$ is simply $k_{n+1}-b_n=\tilde{k}_n-b_n$, making for a total contribution over all depths of $(-R)^{\tilde{k}-b}$.  

Finally, each instance of transmission from the $j$th layer to the $(j+1)$st layer in $\mpp$ corresponds to a vertex in the tree representing $\mpp$ which is not a leaf, i.e., to a branch point at $z_j$; and, since the path $\mpp$ starts and ends at $z_{-1}$ every such transmission has a corresponding return transmission in the opposite direction, from the $(j+1)$st layer to the $j$th layer.  There are $b_j$ branch points at $z_j$, and each of these corresponds to two transmissions across the boundary at $z_j$, making for a total contribution to $w(\mpp)$ of $T^{2b}$.   Since every $w_n$ is covered by one of the above cases, the lemma follows. 
\end{pf}

\subsection{The Green's function}

A formula for the coefficients $a(R,k)$ in (\ref{arrival}), defined as the summation (\ref{aRkdefintion}), is now within easy reach.  Recall that the binomial coefficient $\binom{x}{y}$ for a pair of non-negative integer vectors $x,y\in\integer^{M+1}$, with $y\leq x$, is to be interpreted as 
\[
\binom{x}{y}=\prod_{n=0}^M\binom{x_n}{y_n}.
\]
(The inequality $y\leq x$ means that $x-y$ has non-negative entries.)
\begin{lem}\label{lem-branch}
Let $k\in\mathfrak{L}_M$ and let $\min\{\mathbb{1},\tilde{k}\}\leq b\leq\min\{k,\tilde{k}\}$.  Then
\[
\#\left\{\mpp\in\sss_M\,|\,(\kappa(\mpp),\beta(\mpp))=(k,b)\right\}=\binom{k}{b}\binom{\tilde{k}-u}{b-u}.
\]
\end{lem}
\begin{pf}
To count the number of scattering sequences having a given transit count vector, it suffices to count the number of corresponding trees---which is straightforward.  
Consider first the arrangement of vertices in a tree for which $(\kappa(\mpp),\beta(\mpp))=(k,b)$.   At each depth $z_n$, there are $k_n$ vertices of which $b_n$ are branch points and $k_n-b_n$ are leaves.  There are $\binom{k_n}{b_n}$ ways of arranging these from left to right, making for a total of 
\begin{equation}\label{vertex}
\binom{k}{b}=\prod_{j=0}^M\binom{k_n}{b_n}
\end{equation}
possible vertex arrangements.  (There is only one way to place the root at $z_{-1}$, which may be ignored.)  

For each vertex arrangement there are various possible edge arrangements, as follows.   Each of the $k_{n+1}=\tilde{k}_n$ vertices at $z_{n+1}$ must be connected by an edge to one of the $b_n$ branch points at $z_n$, respecting the vertex ordering (so that edges don't cross).  This is equivalent to choosing a $b_n$-part ordered partition of the integer $\tilde{k}_n$.  If $\tilde{k}_n\geq 1$, there are $\binom{\tilde{k}_n-1}{b_n-1}$ possible choices; and if $\tilde{k}_n=0$ then $b_n=0$ and there is $1=\binom{\tilde{k}_n}{b_n}$ (empty) arrangement.  Letting $N$ denote the largest index for which $\tilde{k}_N\neq 0$, the total number of edge arrangements is 
\begin{equation}\label{edge}
\binom{\tilde{k}-u}{b-u}=\prod_{n=0}^N\binom{\tilde{k}_n-1}{b_n-1}.
\end{equation}
Combining (\ref{vertex}) and (\ref{edge}) yields a total tree count of 
\[
\#\left\{\mpp\in\sss_M\,|\,(\kappa(\mpp),\beta(\mpp))=(k,b)\right\}=\binom{k}{b}\binom{\tilde{k}-u}{b-u},
\]
completing the proof. 
\end{pf}

Combining the foregoing lemmas gives a formula for the reflection Green's function, as follows. 
\begin{thm}\label{thm-reflection} Let $(\tau,R)$ be an $M$-layer medium for some integer $M\geq 1$. Then 
\[
\gtr(t)=\sum_{k\in\mathfrak{L}_M}a(R,k)\delta(t-\langle k,\tau\rangle)
\]
where for each $k\in\mathfrak{L}_M$, the amplitude $a(R,k)$ is given by the following formula. Setting $u=\min\{\mathbb{1},\tilde{k}\}$,
\begin{equation}\label{aRk}
a(R,k)=\sum_{b\in V(k)}\binom{k}{b}\binom{\tilde{k}-u}{b-u}(-R)^{\tilde{k}-b}R^{k-b}T^{2b},
\end{equation}
where $V(k)$ denotes the set of $b\in\integer^{M+1}$ such that $u\leq b\leq\min\{k,\tilde{k}\}$.
\end{thm}
\begin{pf}
The total amplitude resulting from scattering sequences having a given transit count vector $k\in\mathfrak{L}_M$ is 
\[
a(R,k)=\sum_{\substack{\mpp\in\sss_M\; \scriptstyle{\rm such\; that}\\ \kappa(\mpp)=k}}w(\mpp).
\]
By Proposition~\ref{prop-V} and Lemma~\ref{lem-weight} the above sum may be rearranged as 
\begin{eqnarray*}
a(R,k)&=&\sum_{b\in V(k)}\sum_{\substack{\mpp\in\sss_M\; \scriptstyle{\rm such\; that}\\ (\kappa(\mpp),\beta(\mpp))=(k,b)}}w(\mpp)\\
&=&\sum_{b\in V(k)}\Bigl(\#\{\mpp\in\sss_M\,|\,(\kappa(\mpp),\beta(\mpp))=(k,b)\}\Bigr)(-R)^{\tilde{k}-b}R^{k-b}T^{2b}.
\end{eqnarray*}
Applying Lemma~\ref{lem-branch} then gives the stated formula.  
\end{pf}

Note that since each transmission coefficient $T_n=\sqrt{1-R_n^2}$ occurs to an even power in (\ref{aRk}), the amplitude $a(R,k)$ is a polynomial in the variables $R_0,\ldots,R_M$ of precise degree
\[
\langle k+\tilde{k},\mathbb{1}\rangle.
\]
Note further that this polynomial depends only on the $R_n$ corresponding to the support of the vector $k$: if $k_n=0$, then $a(R,k)$ does not depend on $R_n$.  Because the amplitudes are polynomial functions of the reflection coefficients, Theorem~\ref{thm-reflection} is straightforward to code.  This makes it possible to compute $\gtr(t)$ exactly up to some cutoff time $T>0$, and to do so very efficiently---see Section~\ref{sec-example}.  

The simplest possible transit count vectors are those that correspond to primary scattering sequences (\ref{primary}); given $n\geq 1$, write $k^n$ for the primary transit count vector defined as
\[
k^n_j=\left\{\begin{array}{cc} 
1&\mbox{ if }0\leq j\leq n\\
0&\mbox{ if }n<j
\end{array}
\right..
\]
The formula for the corresponding amplitude is easy to work out directly. It is the only explicit example we found in the existing literature, appearing in the early work of Kunetz \cite{Ku:1963}.  Observe that $V(k^n)=\{k^{n-1}\}$ is a singleton, so that the general formula (\ref{aRk}) reduces just to
\[
a(R,k^n)=(-R)^0R^{k^n-k^{n-1}}T^{2k^{n-1}}=R_n\prod_{j=0}^{n-1}(1-R_j^2),
\]
which is exactly Kunetz's formula.

\section{The combinatorics of transmission}
The case of transmission is in many ways similar to that of reflection.   We therefore give a much more compressed presentation, emphasizing only those points where there is a substantial difference.  

\subsection{Scattering sequences and weights}
In the case of transmission, a scattering sequence is a path in the graph
\begin{equation}\label{simplegraph2}
\mbox{
\includegraphics[width=4in]{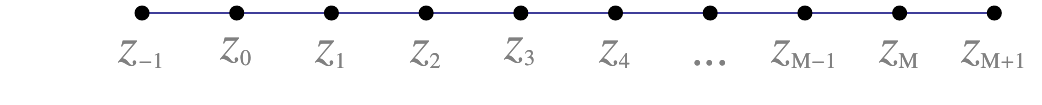}
}
\end{equation}
that starts at $z_{-1}$ and ends at $z_{M+1}$.  (See Figure~\ref{fig-Layered}.)  Formally, given an integer $M\geq 1$, let $\sss_M^\prime$ denote the set of all sequences of the form
\[
\mpp=(\mpp_0,\mpp_1,\ldots,\mpp_L)
\]
such that: 
\begin{subequations}
\begin{align}
\mpp_0=z_{-1},\;\mpp_L=z_{M+1}\mbox{ and }\mpp_n\in&\left\{z_1,\ldots,z_M\right\}\mbox{ if }1\leq n\leq L-1;\label{transpathcondition1}\\
\begin{split}
\forall n \mbox{ with } 0\leq n\leq L-1,&\quad\exists j \mbox{ with }-1\leq j\leq M\mbox{ such that }\\
&\{\mpp_n,\mpp_{n+1}\}=\{z_j,z_{j+1}\}.\label{transpathcondition2}
\end{split}
\end{align}
\end{subequations}
The elements of $\sss_M^\prime$ will be referred to as scattering sequences (or \emph{transmission} scattering sequences for emphasis).  

The weight of a scattering sequence in $\sss_M^\prime$ is computed exactly as for a scattering sequence in $\sss_M$; see Section~\ref{sec-weight}.

\subsection{Trees and transit count vectors}
\begin{figure}[h]
\fbox{
\includegraphics[clip,trim=0in 0in 0in 0in, width=4.8in]{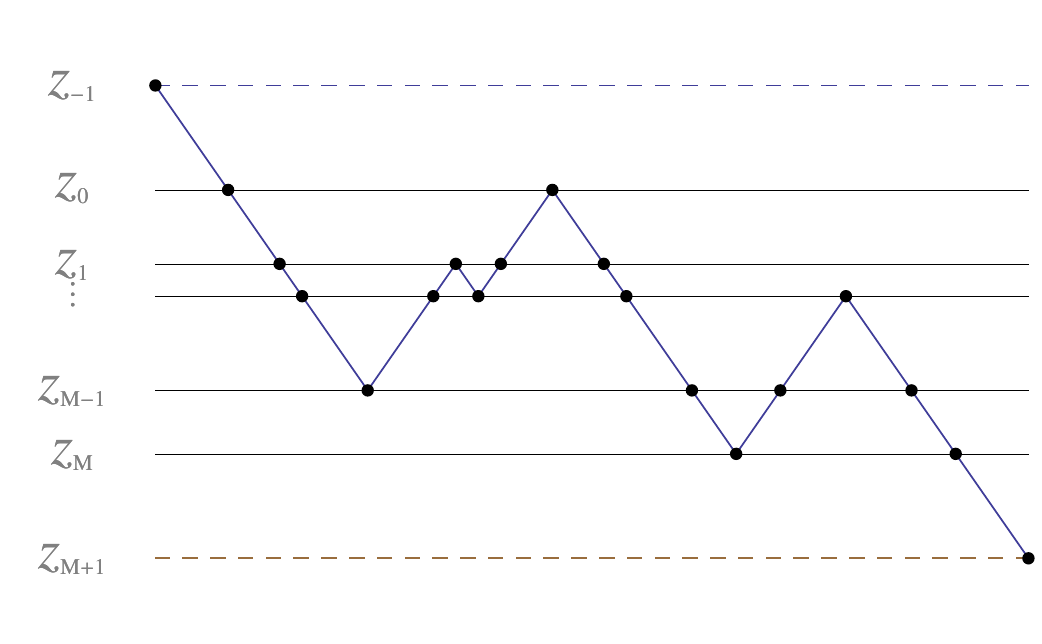}
}
\caption{ The scattering path for a transmission scattering sequence. As with a Dyck path, the horizontal coordinate represents time, increasing to the right.}\label{fig-transmission}
\end{figure}

Figure~\ref{fig-transmission} depicts what we call a \emph{scattering} path, the analogue of a Dyck path for a transmission scattering sequence.  

Note that there are subpaths of a transmission scattering path that themselves are Dyck paths. We call these Dyck subpaths.   
\begin{figure}[h]
\fbox{
\includegraphics[clip,trim=0in 0in 0in 0in, width=2.3in]{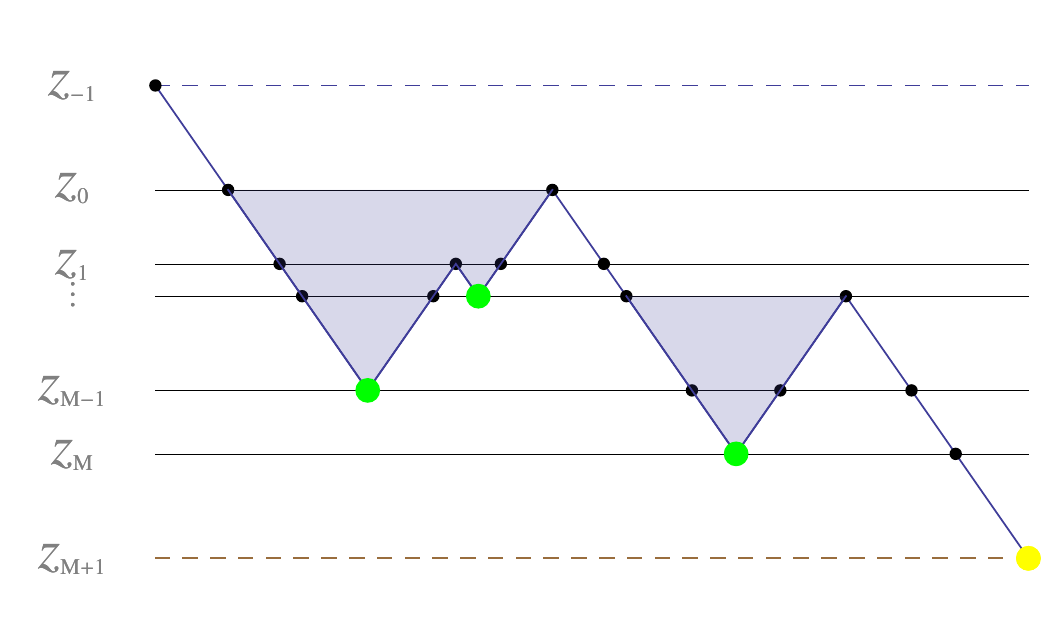}
}\ \fbox{
\includegraphics[clip,trim=0in 0in 0in 0in, width=2.3in]{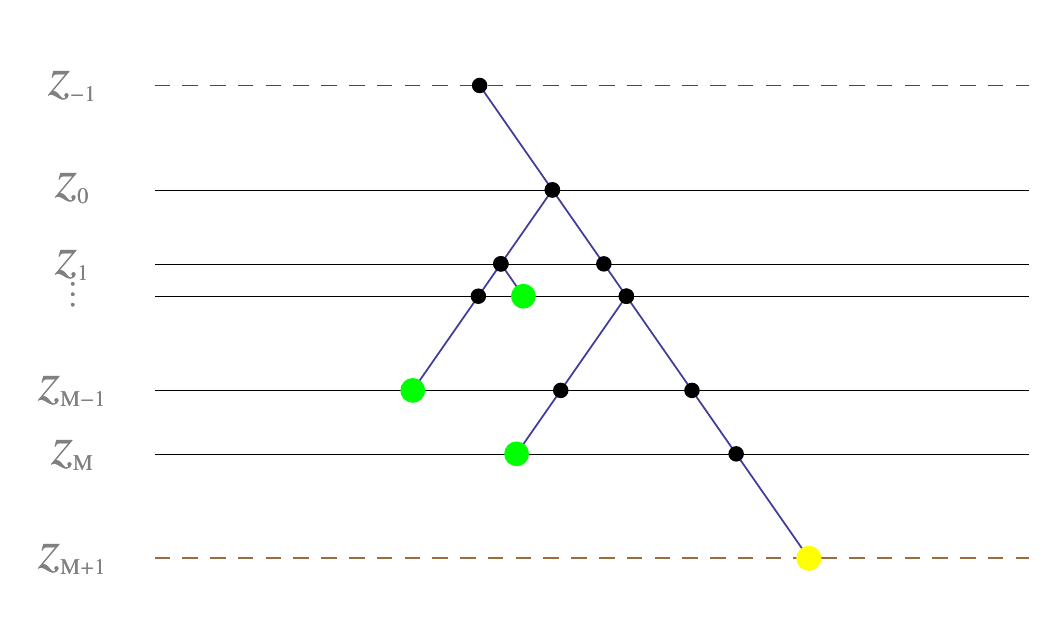}
}

\caption{ The scattering path with Dyck subpaths filled in, and it's associated tree (depicted upside down, with the root at the top) obtained by collapsing the Dyck subpaths.  The tip is coloured yellow; the other leaves are green.}\label{fig-subpaths}
\end{figure}
Collapsing the Dyck subpaths turns the scattering path into a tree, which has some additional structure.   There is a distinguished leaf, the single leaf of maximum height (corresponding to $z_{M+1}$), which we call the \emph{tip}.  The subpath leading directly from the root to the tip will be referred to as the \emph{trunk}.   See Figure~\ref{fig-subpaths}.

Collapsing of Dyck paths is reversible, so there is a one-to-one correspondence between scattering paths and trees; the latter serve as a convenient representation for the purpose of counting.   Since every tree that corresponds to a scattering sequence has a trunk (including the root and the tip of the tree), it simplifies the combinatorial analysis to focus only on the remaining part of the tree.  As in Section~\ref{sec-reflection}, we define two maps, 
\[
\kappa^\prime,\beta^\prime:\sss_M^\prime\rightarrow\integer^{M+1}
\]
that associate integer vectors to a given transmission scattering sequence $\mpp\in\sss_M^\prime$.   For each $n$ in the range $0\leq n\leq M$, define $k_n$ to be the number of nodes at height $z_n$ not on the trunk.  (Since there is precisely one node at each level belonging to the trunk, $k_n$ is one less than the total number of nodes at height $z_n$.)  Set 
\[
\kappa^\prime(\mpp)=(k_0,k_1,\ldots,k_M).
\]
Note that with this defintion $k_0=0$; the other $k_n$ can take any non-negative integer value.    Define $m_n$ to be the number of branch points at $z_n$ not on the trunk, and set
\[
\beta^\prime(\mpp)=(m_0,\ldots,m_M).
\]
As before, we call $k=(k_0,\ldots,k_M)$ a transit count vector and $m=(m_0,\ldots,m_M)$ a branch count vector, and we define $\tilde{k}$ to be the left shift of $k$:
\[
\tilde{k}=(k_1,k_2,\ldots,k_M,0).
\]
Observe that for each $n$, 
\begin{equation}\label{transmission-constraints}
0\leq m_n\leq\mbox{min}\{k_n,\tilde{k}_n\}.
\end{equation}
Furthermore, letting $\integer_+$ denote the nonnegative integers, the following result is straightforward (we write $0$ for both the number and the zero vector).   
\begin{prop}\label{prop-trans}
Given any $k\in\{0\}\times\integer^M_+$ and any $0\leq m\leq \mbox{min}\{k,\tilde{k}\}$,  there exists a transmission scattering sequence $p\in\sss_M^\prime$ such that 
\[
\bigl(\kappa^\prime(\mpp),\beta^\prime(\mpp)\bigr)=(k,m).
\]
\end{prop}

The following is the analogue of Lemma~\ref{lem-weight} for transmission scattering sequences.
\begin{lem}\label{lem-transmission-weight}
Let $\mpp\in\sss_M^\prime$ be a transmission scattering sequence in an $M$-layer medium $(\tau^\prime,R)$, and set $(k,m)=(\kappa^\prime(\mpp),\beta^\prime(\mpp))$.  Then
\[
w(\mpp)=(-R)^{\tilde{k}-m}R^{k-m}T^{2m+\mathbb{1}}.
\]
\end{lem}
\begin{pf}
Consider the tree representing $\mpp$.  Observe that each instance of $w_n=R_j$ in (\ref{wn}) corresponds to a unique leaf at $z_j$ (see the right-hand part of Figure~\ref{fig-subpaths}).  Since there are $k_n-m_n$ leaves at $z_n$ this results in a total contribution of $R^{k-m}$.   

Let $v$ be a branch point at $z_j$ (possibly on the trunk) having $d_v$ edges to vertices at $z_{j+1}$.   Observe that precisely $d_v-1$ of these edges correspond to an instance of $w_n=-R_j$ in (\ref{wn}), and every occurrence of $w_n=-R_j$ arises this way.  The sum total of numbers $d_v-1$ over branch points $v$ at $z_j$ is simply $(k_{n+1}+1)-(m_n+1)=\tilde{k}_n-m_n$, making for a total contribution over all depths of $(-R)^{\tilde{k}-m}$.  

Finally, each instance of transmission from the $j$th layer to the $(j+1)$st layer in $\mpp$ corresponds to a vertex in the tree representing $\mpp$ which is not a leaf, i.e., to a branch point at $z_j$.  Each branch point not on the trunk corresponds to two transmissions across the boundary at $z_j$ (as for a Dyck path), making for a total contribution to $w(\mpp)$ of $T^{2m}$.  Each branch point on the trunk corresponds to a single upward transmission, making for a total contribution of $T^{\mathbb{1}}$. Since every $w_n$ is covered by one of the above cases, the lemma follows. 
\end{pf}

\subsection{The transmission Green's function}

The analogue of Lemma~\ref{lem-branch} for transmission scattering sequences is slightly simpler, as follows. 
\begin{lem}\label{lem-transmission-branch}
Let $k\in\{0\}\times\integer_+^M$ and let $0\leq m\leq\min\{k,\tilde{k}\}$.  Then
\[
\#\left\{\mpp\in\sss_M^\prime\,|\,(\kappa^\prime(\mpp),\beta^\prime(\mpp))=(k,m)\right\}=\binom{k}{m}\binom{\tilde{k}}{m}.
\]
\end{lem}
\begin{pf}
To count the number of scattering sequences having a given transit count vector, it suffices to count the number of corresponding trees, as before.  
Consider first the arrangement of vertices in a tree for which $(\kappa^\prime(\mpp),\beta^\prime(\mpp))=(k,m)$.  To begin with, there is the trunk, which has fixed structure.   At each depth $z_n$, there are $k_n$ vertices off the trunk, of which $m_n$ are branch points and $k_n-m_n$ are leaves.  There are $\binom{k_n}{m_n}$ ways of arranging these from left to right, making for a total of 
\begin{equation}\label{transmission-vertex}
\binom{k}{m}=\prod_{j=0}^M\binom{k_n}{m_n}
\end{equation}
possible vertex arrangements.  

For each vertex arrangement there are various possible edge arrangements, as follows.   Each of the $k_{n+1}=\tilde{k}_n$ vertices at $z_{n+1}$ must be connected by an edge to one of the $m_n+1$ branch points at $z_n$ (including the branch point on the trunk), respecting the vertex ordering (so that edges don't cross).  This is equivalent to choosing a $m_n+1$-part ordered partition of the integer $\tilde{k}_n$, where the last part---corresponding to the trunk branch point---can be empty.  There are $\binom{\tilde{k}_n}{m_n+1-1}$ possible choices.  Letting $N$ denote the largest index for which $\tilde{k}_N\neq 0$, the total number of edge arrangements is 
\begin{equation}\label{transmission-edge}
\binom{\tilde{k}}{m}=\prod_{n=0}^N\binom{\tilde{k}_n}{m_n}.
\end{equation}
Combining (\ref{transmission-vertex}) and (\ref{transmission-edge}) yields a total tree count of 
\[
\#\left\{\mpp\in\sss_M\,|\,(\kappa(\mpp),\beta(\mpp))=(k,m)\right\}=\binom{k}{m}\binom{\tilde{k}}{m},
\]
as was to be shown.  
\end{pf}

\begin{thm}\label{thm-transmission}
Let $(\tau^\prime,R)$ correspond to an $M$-layer medium for some integer $M\geq 1$, and write $|\tau^\prime|=\tau_0+\cdots+\tau_{M+1}$. Then 
\[
H^{(\tau^\prime,R)}(t)=\sum_{k\in\{0\}\times\integer^M_+}b(R,k)\,\delta\bigl(t-{\scriptstyle\frac{1}{2}}|\tau^\prime|-\langle k,\tau\rangle\bigr)
\]
where for each $k\in\{0\}\times\integer^M_+$, the amplitude $b(R,k)$ is given by the formula\begin{equation}\label{bRk}
b(R,k)=\sum_{0\leq m\leq\min\{k,\tilde{k}\}}\binom{k}{m}\binom{\tilde{k}}{m}(-R)^{\tilde{k}-m}R^{k-m}T^{2m+\mathbb{1}}.
\end{equation}
\end{thm}
\begin{pf}
Given $\tau^\prime$, observe that the arrival time for a transmission scattering sequence $\mpp\in\sss_M^\prime$ having $\kappa^\prime(\mpp)=k$ is precisely 
\[
{\scriptstyle\frac{1}{2}}|\tau^\prime|+\langle k,\tau\rangle.
\]
The theorem then follows from Lemmas~\ref{lem-transmission-weight} and \ref{lem-transmission-branch}.\end{pf}

Note that the factor $T^\mathbb{1}$ occurring in the formula for the transmission amplitude $b(R,k)$ has the form 
\[
T^\mathbb{1}=\sqrt{1-R_0^2}\sqrt{1-R_1^2}\cdots\sqrt{1-R_M^2},
\]
which is not a polynomial in the reflection coefficients $(R_0,\ldots,R_M)$.  So it turns out that $b(R,k)/T^{\mathbb{1}}$ is a polynomial in the reflection coefficients, while $b(R,k)$ itself is not.   (Recall that by contrast reflection amplitudes $a(R,k)$ are polynomial functions of the reflection coefficients.)

\section{An example\label{sec-example}}
\begin{figure}[h]
\fbox{\parbox{\textwidth}{
(a)\includegraphics[clip,trim=0in 0in 0in 0in, width=2.25in]{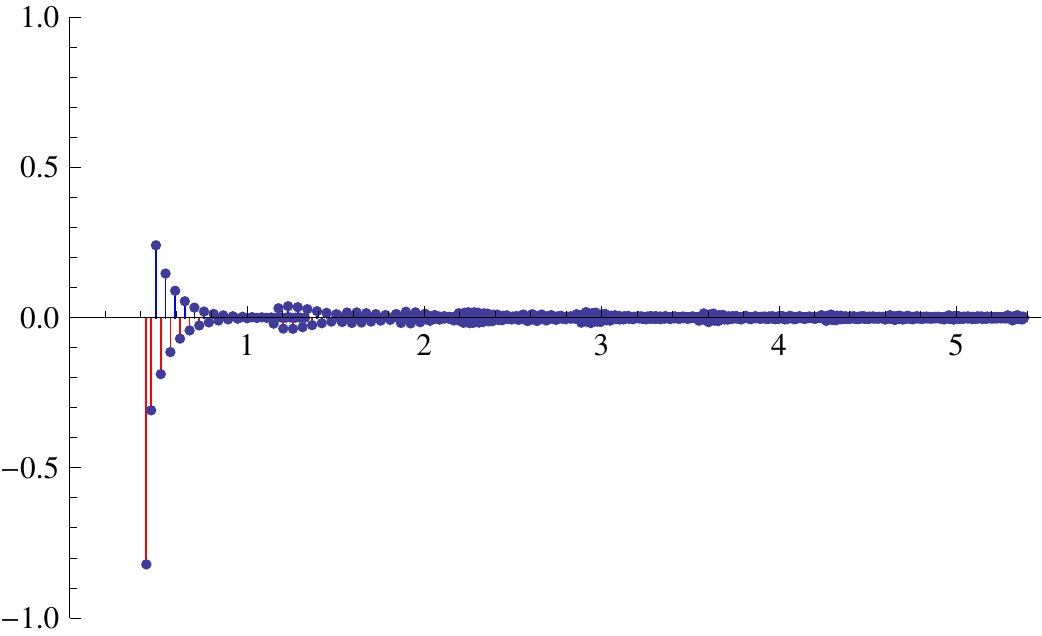}
(b)\includegraphics[clip,trim=0in 0in 0in 0in, width=2.25in]{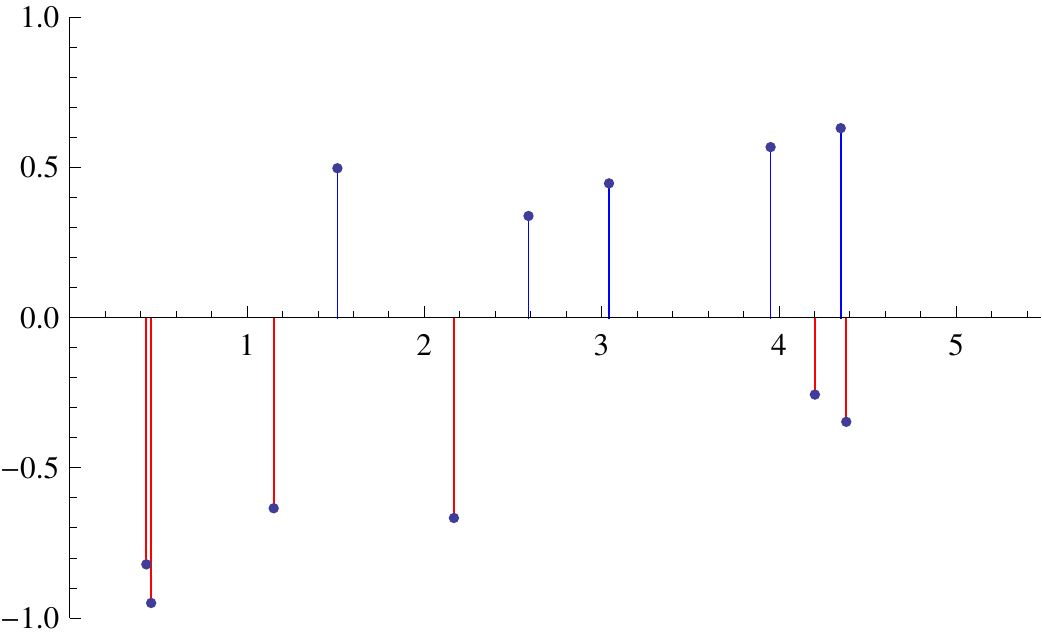}
(c)\includegraphics[clip,trim=0in 0in 0in 0in, width=2.25in]{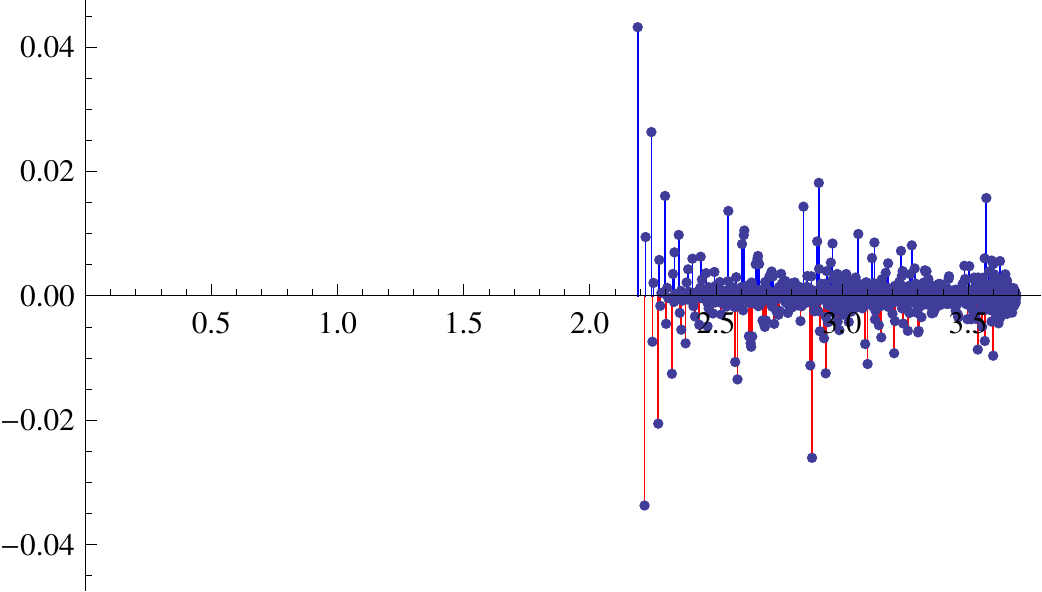}
(d)\includegraphics[clip,trim=0in 0in 0in 0in, width=2.25in]{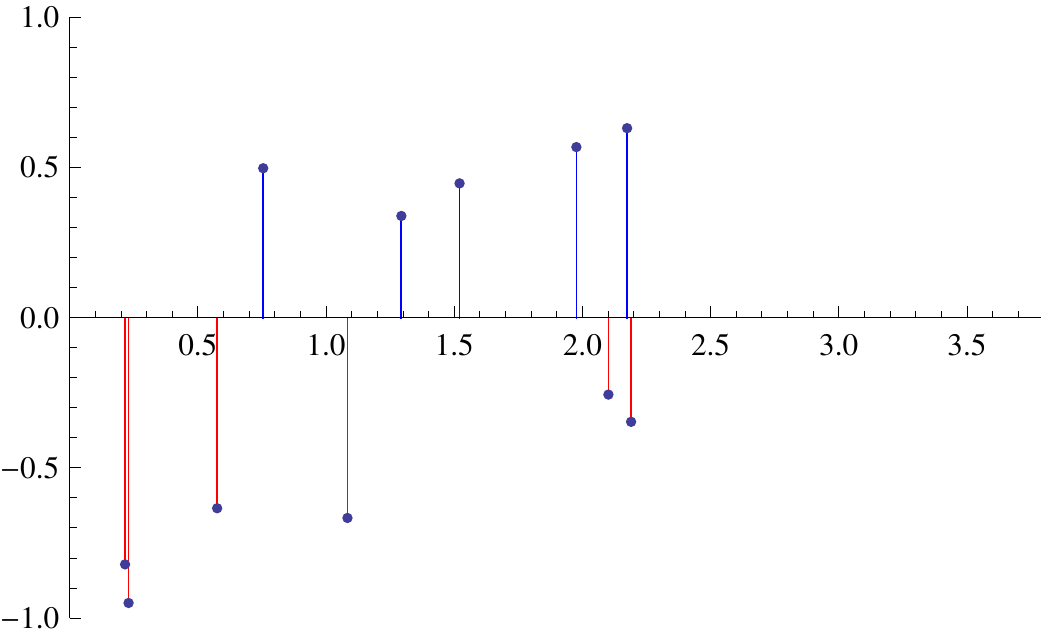}
}
}
\caption{ The exact reflection Green's function (a) for the example (\ref{example}), and the underlying reflection coefficients plotted against two-way travel time (b).  The exact transmission Green's function (c) is plotted on a smaller amplitude scale, with the same reflection coefficients plotted against one-way travel time (d).}\label{fig-example}
\end{figure}
The formulas presented in Theorems~\ref{thm-reflection} and \ref{thm-transmission} make it easy to compute the reflection and transmission Green's functions exactly up to a finite cutoff time $T>0$.   We illustrate this in the present situation by working out a numerical example.  Consider the pair $(\tau,R)$ representing an $M=10$ layer medium, with the following values.    
\begin{equation}\label{example}
\begin{array}{cccc}
n&\mbox{\parbox{1.45in}{Two-way travel time from $z_{-1}$ to $z_n$}}& \tau_n &R_n\\
\hline \\
0& 0.432779 & 0.432779 & -0.821708 \\
 1&0.459894 & 0.0271153 & -0.950247 \\
2& 1.1509 & 0.69101 & -0.634818 \\
3& 1.51027 & 0.35937 & 0.497529 \\
 4&2.16748 & 0.657201 & -0.66719 \\
 5&2.58797 & 0.420491 & 0.338592 \\
 6&3.04205 & 0.454083 & 0.447163 \\
 7&3.95335 & 0.911298 & 0.567927 \\
 8&4.20316 & 0.249813 & -0.256277 \\
 9&4.34889 & 0.145731 & 0.630965 \\
 10&4.38014 & 0.0312502 & -0.346928
\end{array}
\end{equation}
Setting $\tau_{M+1}=0$, which corresponds to measuring the transmission response exactly at $z_M$, formulas for $G^{(\tau,R)}$ and $H^{(\tau^\prime,R)}$ given in Theorems~\ref{thm-reflection} and \ref{thm-transmission} were coded up in Mathematica.  The resulting functions were computed up to a cutoff time of 5.38 seconds for reflection and 3.69 seconds for transmission.  (The respective computations took 3.37 seconds and 17.35 seconds on a 1.8 GHz i7 processor.) The computed pulse trains are depicted in Figure~\ref{fig-example}.   The reflected pulse train has 19~242 terms, and the transmitted one has 35~059 terms;  the majority of these are of an amplitude which is too small to be visible in the plots.   

This example illustrates that Theorems~\ref{thm-reflection} and \ref{thm-transmission},  while derived by combinatorial methods not usually associated with the analysis of PDEs, offer a computationally tractable representation of reflected and transmitted waves.   This provides a new tool for application to the various imaging modalities where one-dimensional models are important.

\appendix

\section{Appendix: The underlying PDE\label{sec-appendix}}

Let $u(t,z)$ denote the velocity (in the $z$-direction) of a material particle at depth $z$ and time $t$, and let $p(t,z)$ denote the pressure.   The medium evolves according to the coupled one-dimensional equations
\begin{subequations}
\begin{align}
\rho\frac{\partial u}{\partial t}+\frac{\partial p}{\partial z}&=0\label{wave1}\\
\frac{1}{K}\frac{\partial p}{\partial t}+\frac{\partial u}{\partial z}&=0.\label{wave2}
\end{align}
\end{subequations}
where $\rho(z)$ is the material density at depth $z$ and $K(z)$ is the bulk modulus.  
For the sake of definiteness we focus on the velocity field $u(t,z)$, although the results can just as easily be formulated in terms of $p(t,z)$.
The initial conditions corresponding to a plane wave unit impulse propagating downward from $z_{-1}$ are
\begin{equation}\label{initial}
\begin{split}
u(0,z)&=\delta(z-z_{-1})\\
 p(0,z)&=\sqrt{K(z_{-1})\rho(z_{-1})}\;\delta(z-z_{-1}).
 \end{split}
\end{equation}
Letting $G(t)$ denote the (velocity) impulse response at $z_{-1}$, we have that  
\begin{equation}\label{response}
G(t)=u(t,z_{-1})\quad\quad(t>0),
\end{equation}
the solution at depth $z_{-1}$ to the system (\ref{wave1},\ref{wave2},\ref{initial}).   Similarly, the impulse response at $z_{M+1}$ is 
\begin{equation}\label{transmission-response}
H(t)=u(t,z_{M+1})\quad\quad(t>0),
\end{equation}
the solution at depth $z_{M+1}$ to the system (\ref{wave1},\ref{wave2},\ref{initial}).

The following facts are derived in  \cite[Chapter~3]{FoGaPaSo:2007}, \cite[Section~2]{BuBu:1983} and elsewhere.   For $1\leq n\leq M$, let $\tau_n$ denote the two-way travel time (for a traveling wave) across the $n$th layer of the above $M$-layer medium, and let $\tau_0$ denote the two-way travel time from depth $z_{-1}$ to $z_0$.   For $0\leq n\leq M$, let $R_n$ denote the reflection coefficient at depth $z_n$ relative to a wave traveling toward the interface from above.  Letting $K_n$ and $\rho_n$ denote the density and bulk modulus inside the $n$th layer---with $K_{-1},\rho_{-1}$ and $K_{M+1},\rho_{M+1}$ 
denoting the respective values at $z_{-1}$ and any point $z_{M+1}$ below $z_{M}$---the travel times and reflection coefficients are given by the formulas
\begin{equation}\label{tau-R-formulas}
\tau_n=\frac{2(z_n-z_{n-1})}{\sqrt{K_n/\rho_n}}\quad\mbox{ and }\quad R_n=\frac{\sqrt{K_n\rho_n}-\sqrt{K_{n+1}\rho_{n+1}}}{\sqrt{K_n\rho_n}+\sqrt{K_{n+1}\rho_{n+1}}},
\end{equation}
for $0\leq n\leq M$.  Note that $-1<R_n<1$ by virtue of (\ref{tau-R-formulas}).

\bibliography{ReferencesEAL}

\begin{thebibliography}{10}

\bibitem{BlCoSt:2001}
N.~Bleistein, J.~K. Cohen, and J.~W. Stockwell, Jr.
\newblock {\em Mathematics of multidimensional seismic imaging, migration, and
  inversion}, volume~13 of {\em Interdisciplinary Applied Mathematics}.
\newblock Springer-Verlag, New York, 2001.
\newblock Geophysics and Planetary Sciences.

\bibitem{BrGo:1990}
L.~M. Brekhovskikh and O.~A. Godin.
\newblock {\em Acoustics of Layered Media {I}}, volume~5 of {\em Springer
  Series on Wave Phenomena}.
\newblock Springer, Heidelberg, 1990.

\bibitem{BuBu:1983}
K.~P. Bube and R.~Burridge.
\newblock The one-dimensional inverse problem of reflection seismology.
\newblock {\em SIAM Rev.}, 25(4):497--559, 1983.

\bibitem{FoGaPaSo:2007}
J.-P. Fouque, J.~Garnier, G.~Papanicolaou, and K.~S{\o}lna.
\newblock {\em Wave propagation and time reversal in randomly layered media},
  volume~56 of {\em Stochastic Modelling and Applied Probability}.
\newblock Springer, New York, 2007.

\bibitem{In:2009}
K.~A. Innanen.
\newblock Born series forward modelling of seismic primary and multiple
  reflections: an inverse scattering shortcut.
\newblock {\em Geophysical Journal International}, 177(3):1197--1204, 2009.

\bibitem{Ku:1963}
G.~Kunetz.
\newblock Quelques exemples d'analyse d'enregistrements sismiques.
\newblock {\em Geophysical Prospecting}, 11(4):409--422, 1963.

\bibitem{Le:2005}
J.-F. Le~Gall.
\newblock Random trees and applications.
\newblock {\em Probab. Surv.}, 2:245--311, 2005.

\bibitem{Ne:1980}
R.~G. Newton.
\newblock Inversion of reflection data for layered media: a review of exact
  methods.
\newblock {\em Geophysical Journal of the Royal Astronomical Society},
  65(1):191--215, 1981.

\bibitem{SaSy:1988}
F.~Santosa and W.~W. Symes.
\newblock Reconstruction of blocky impedance profiles from normal-incidence
  reflection seismograms which are band-limited and miscalibrated.
\newblock {\em Wave Motion}, 10(3):209--230, 1988.

\bibitem{St:1999}
R.~P. Stanley.
\newblock {\em Enumerative combinatorics. {V}ol. 2}, volume~62 of {\em
  Cambridge Studies in Advanced Mathematics}.
\newblock Cambridge University Press, Cambridge, 1999.
\newblock With a foreword by Gian-Carlo Rota and appendix 1 by Sergey Fomin.

\bibitem{WeArFe:2003}
A.~B. Weglein, F.~V. Ara{\'u}jo, P.~M. Carvalho, R.~H. Stolt, K.~H. Matson,
  R.~T. Coates, D.~Corrigan, D.~J. Foster, S.~A. Shaw, and H.~Zhang.
\newblock Inverse scattering series and seismic exploration.
\newblock {\em Inverse Problems}, 19(6):R27--R83, 2003.

\end{thebibliography}

\end{document}